# Nonparametric estimation of a convex bathtub-shaped hazard function

HANNA K. JANKOWSKI[1] and JON A. WELLNER[2]

[1]*Department of Mathematics and Statistics, N520 Ross Building, 4700 Keele Street, York University, Toronto, ON, Canada M3J 1P3.*
*E-mail: hkj@mathstat.yorku.ca*

[2]*Department of Statistics, Box 354322, University of Washington, Seattle, WA 98195-4322 USA.*
*E-mail: jaw@stat.washington.edu*

In this paper, we study the nonparametric maximum likelihood estimator (MLE) of a convex hazard function. We show that the MLE is consistent and converges at a local rate of $n^{2/5}$ at points $x_0$ where the true hazard function is positive and strictly convex. Moreover, we establish the pointwise asymptotic distribution theory of our estimator under these same assumptions. One notable feature of the nonparametric MLE studied here is that no arbitrary choice of tuning parameter (or complicated data-adaptive selection of the tuning parameter) is required.

## 1. Introduction

Information on the behavior of time to a random event is of much interest in many fields. The random event could be failure of a material or machine, death, an earthquake or infection by a disease, to name but a few examples. Frequently, this type of data is called *lifetime data*, and it is natural to assume that it takes values in $[0, \infty)$. If the lifetime distribution $F$ has a density $f$, then a key quantity of interest is the hazard (or failure) rate $h(t) = f(t)/(1 - F(t))$. Heuristically, $h(t) \, \mathrm{d}t$ is the probability that, given survival until time $t$, the event will occur in the next duration of length $\mathrm{d}t$. The hazard function is also known as the *force of mortality* in actuarial science or the *intensity function* in extreme value theory.

Certain shape restrictions arise quite naturally for hazard rates. In this work, we are particularly interested in the family of hazard functions which are convex. That is, we attach the additional smoothness constraint of convexity to the more traditional assumption of a bathtub-shaped failure rate (that is, first decreasing, then increasing). Heuristically, bathtub-shaped hazards correspond to lifetime distributions with high initial hazard (or infant mortality), lower and often rather constant hazard during the middle of life and then increasing hazard of failure (or wear-out) as aging proceeds; see [20, 23].

Many other estimators of hazard functions (and solutions to the closely related problem of estimating the intensity of a Poisson process) with and without shape restrictions have been considered in the literature; see [24] for a partial review up to 2002. In recent years, the focus has shifted to construction of "adaptive" estimators over large scales of smoothness classes; see, for example, [5, 6, 25]. Virtually all of these other estimators require careful choice of penalty terms or tuning parameters, and computation of the adaptive estimators typically involves methods of combinatorial optimization. As far as we know, reliable algorithms for computing them are not yet available. Our estimators avoid the choices of tuning parameters or penalty terms by virtue of the shape constraint of convexity and are relatively straightforward to compute since the corresponding optimization problems are convex.

Recall the definition of a convex function. Let $C \subset \mathbb{R}_+ = [0, \infty)$ be convex. Then $h : C \mapsto \mathbb{R}$ is convex (on $C$) if it satisfies

$$h(\lambda x + (1-\lambda)y) \leq \lambda h(x) + (1-\lambda)h(y), \qquad 0 < \lambda < 1,$$

for all $x, y \in C$. Equivalently, a function is convex if its *epigraph*

$$\{(x, \mu) : x \in C, \mu \in \mathbb{R}, \mu \geq f(x)\}$$

is a convex set in $\mathbb{R}^2$ (see, for example, [27], Section 4). Thus, a convex function on $C$ may be extended to a convex function on $\mathbb{R}_+$ by setting $h(x) = +\infty$ for $x \in \mathbb{R}_+ \cap C^c$.

Suppose, then, that we observe i.i.d. variables $X_1, \ldots, X_n$ from a distribution $F_0$ with density $f_0$ and hazard rate $h_0$. We denote the true cumulative hazard function by $H_0(t) = \int_0^t h_0(s) \, ds$, and the true survival function by $S_0 = 1 - F_0$. Also, $0 < X_{(1)} < X_{(2)} < \cdots < X_{(n)}$ denote the order statistics corresponding to $X_1, \ldots, X_n$.

To define the MLE of $h_0$, $\widehat{h}_n$, we first consider the likelihood in terms of the hazard,

$$\mathcal{L}(h) = \prod_{i=1}^n h(X_i) \exp\{-H(X_i)\} = \prod_{i=1}^n h(X_{(i)}) \exp\{-H(X_{(i)})\},$$

where $H(t) = \int_0^t h(s) \, ds$. This can be made arbitrarily large by increasing the value of $h(X_{(n)})$. We therefore find $\widehat{h}_n : [0, X_{(n)}) \mapsto \mathbb{R}_+$ by maximizing the modified likelihood

$$\mathcal{L}^{\mathrm{mod}}(h) = \prod_{i=1}^{n-1} h(X_{(i)}) \exp\{-H(X_{(i)})\} \times \exp\{-H(X_{(n)})\} \qquad (1.1)$$

over $\mathcal{K}_+$, the space of non-negative convex functions on $[0, X_{(n)})$. The *full* MLE is then found by setting $\widehat{h}_n(x) = \infty$ for all $x \geq X_{(n)}$. This is the same approach as taken in [10]. Equivalently, one could first impose the constraint that $h \leq M$, and then let $M \to \infty$ (see, for example, [26], page 338).

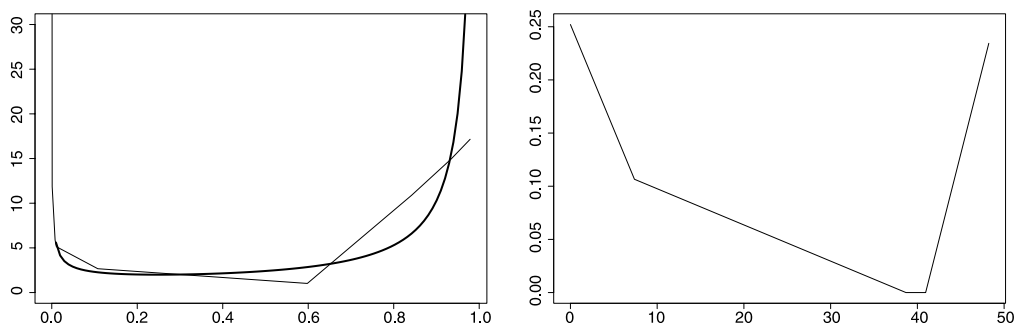

**Figure 1.** Examples of the estimator. Left: Estimating the HS hazard with $b=0$ and $A=1$ for a sample size of 100 (bold = true hazard, solid = MLE). Right: Estimation of the earthquake hazard from CPTI04 data (solid = MLE).

To illustrate the proposed estimator, consider the distribution with density given by

$$f(t) = \frac{1+2b}{2A\sqrt{b^2 + (1+2b)t/A}} \qquad \text{on } 0 \le t \le A.$$

This distribution was derived in [14] as a relatively simple model with bathtub-shaped hazards which also has an adequate ability to model lifetime behavior. We will call this the *HS distribution*, after the authors. It has convex hazards for all values of $b$ in the parameter space ($b > -1/2$). Figure 1 shows the MLE for simulated data from this distribution with sample size $n = 100$.

We also applied our estimators to the earthquake data of the *Appennino Abruzzese* region of Italy (Region 923) recently considered by [19], where Bayesian estimation methods are studied. The data comes from the Gruppo di Lavoro CPTI (2004) catalog [9]. It consists of 46 inter-quake times for Region 923, occurring after the year 1650 and with moment magnitude greater than 5.1 (details on the justification of these criteria is available in [19], page 14). Figure 1 shows the resulting estimator.

The main results of this paper are the characterizations, consistency and asymptotic behavior of the nonparametric MLE of a convex hazard function. The estimator is continuous and piecewise linear on $[0, X_{(n)})$. Although we give a characterization of the MLE, the final form of the estimator is not explicit. We therefore propose an algorithm (based on the support reduction algorithm of [13]). This algorithm is discussed in a separate report, [16], and is available as the R package convexHaz [17].

To describe the local asymptotics of the MLE, we introduce the following process.

**Definition 1.1.** *Let $W(s)$ denote a standard two-sided Brownian motion, with $W(0) = 0$, and define $Y(t) = \int_0^t W(s)\,\mathrm{d}s + t^4$. The function $\{\mathcal{I}(t): t \in \mathbb{R}\}$, the invelope of the process $\{Y(t): t \in \mathbb{R}\}$, is defined as follows:*

$$\text{the function } \mathcal{I} \text{ is above the function } Y: \mathcal{I}(t) \ge Y(t) \text{ for all } t \in \mathbb{R}; \qquad (1.2)$$

*the function $\mathcal{I}$ has a convex second derivative;* (1.3)

*the function $\mathcal{I}$ satisfies $\int_{\mathbb{R}} \{\mathcal{I}(t) - Y(t)\} \, d\mathcal{I}^{(3)}(t) = 0$.* (1.4)

It was shown in [11] that the process $\mathcal{I}$ exists and is almost surely uniquely defined. Moreover, with probability one, $\mathcal{I}$ is three times differentiable at $t=0$. The asymptotic behavior of all of our estimators may be described in terms of the derivatives of the envelope $\mathcal{I}$ at zero. The following theorem builds on the basic results in [12] concerning nonparametric estimation of a decreasing convex density.

**Theorem 1.2.** *Suppose that $h_0$ is convex, $x_0 > 0$ is a point which satisfies $h_0(x_0) > 0$ and $h_0''(x_0) > 0$, and $h_0''(\cdot)$ is continuous in a neighborhood of $x_0$. Then*

$$\begin{pmatrix} n^{2/5}(\widehat{h}_n(x_0) - h_0(x_0)) \\ n^{1/5}(\widehat{h}_n'(x_0) - h_0'(x_0)) \end{pmatrix} \to_d \begin{pmatrix} c_1 \mathcal{I}^{(2)}(0) \\ c_2 \mathcal{I}^{(3)}(0) \end{pmatrix},$$

*where $\mathcal{I}^{(2)}(0)$ and $\mathcal{I}^{(3)}(0)$ are the second and third derivatives at 0, respectively, of the envelope of $Y(t) \equiv \int_0^t W(s) \, ds + t^4$ and where*

$$c_1 = \left( \frac{h_0^2(x_0) h_0''(x_0)}{24 S_0^2(x_0)} \right)^{1/5}, \qquad c_2 = \left( \frac{h_0(x_0) h_0''(x_0)^3}{24^3 S_0(x_0)} \right)^{1/5}.$$

The key to this result lies in Lemma 5.3, where we establish that the "touchpoints" (defined carefully in Section 2) cluster around $x_0$ at a local scale of $n^{-1/5}$. The assumption that $h_0''$ is strictly convex and continuous near $x_0$ is crucial in this step. If $h_0''(x_0) = 0$ (and is continuous in a neighborhood of $x_0$), we conjecture that $\widehat{h}_n(x_0)$ converges at the rate $n^{1/2}$. Similar behavior has been noted for monotone density estimators in [8]. If $h_0'$ is discontinuous at $x_0$, unpublished work of Cai and Low [7] suggests that $\widehat{h}_n(x_0)$ converges to $h_0(x_0)$ at rate $n^{1/3}$. The behavior of convex-constrained estimators in both of these situations remains unknown and is the subject of current research.

The limiting distributions of $\widehat{h}_n(x_0)$ and $\widehat{h}_n'(x_0)$ involve the constants $c_1$ and $c_2$, which depend on the (unknown) hazard function $h_0$, as well as the random variables $(\mathcal{I}^{(2)}(0), \mathcal{I}^{(3)}(0))$ which have a universal distribution free of the parameters of the problem. Thus, Theorem 1.2 can be used, in principle, to form confidence intervals for $h_0(x_0)$ and $h_0'(x_0)$. This would involve estimation of the constants $c_1 = c_1(h_0, x_0)$ and $c_2 = c_2(h_0, x_0)$, respectively, both of which depend on $h_0''(x_0)$, and appropriate quantiles of the distributions of $\mathcal{I}^{(2)}(0)$ and $\mathcal{I}^{(3)}(0)$, respectively. Although virtually nothing is known about the distribution of the invelope and its derivatives analytically, the algorithms developed in [11] can easily be used to obtain simulated values of the needed quantiles. Other possible approaches to confidence intervals in this problem involve inversion of likelihood ratio tests (see [2–4] for this approach in the context of monotone or $U$-shaped function estimation) or resampling methods as in [22] and as discussed in [4] in the setting of nonparametric estimation of monotone functions. It should be noted that our Theorem 1.2 verifies one of the key hypotheses needed for validity of the general subsampling theory of [21, 22], and therefore makes the subsampling approach to confidence intervals viable. The details and properties of all these approaches remain to be investigated.

The outline of this paper is as follows. Section 2 is dedicated to the proof of characterizations, existence and uniqueness of the MLE. Consistency is proved in Sections 3 and 4 establishes lower bounds for the pointwise minimax risk of $\widehat{h}_n$. Rates of convergence are established in Section 5, with Section 5.3 containing proofs of our main results concerning the limiting distribution at a fixed point. The companion technical report [15] also includes a detailed treatment of a least-squares estimator, as well as sketches of similar results for censored data and intensity functions of Poisson processes.

## 2. Characterizations, uniqueness and existence

**Proposition 2.1.** *The function $\widehat{h}_n$ which minimizes $\varphi_n$ over $\mathcal{K}_+$ is piecewise linear. It has at most one change of slope between observations, except perhaps in one such interval, where, if the estimator touches zero, it may have two changes of slope (it is zero between these two changes). Also, between zero and $X_{(1)}$, the minimizer may have at most one change of slope, but this happens only if it touches zero, and in this case the estimator is increasing and equal to zero before the first change of slope. Between $X_{(n-1)}$ and $X_{(n)}$, the minimizer will also have at most one change of slope, and this only in the case where it is decreasing on $[X_{(n-1)}, X_{(n)})$ and equal to zero after the change.*

**Proof.** Consider any $h$ and choose a convex $g$ such that $h(X_i) = g(X_i)$ for $i = 1, \ldots, n-1$ and $h(x) \geq g(x) \geq 0$ on $[0, X_{(n)})$. It follows that $\varphi_n(h) - \varphi_n(g) \geq 0$ if and only if $H(X_i) \geq G(X_i)$ for $i = 1, \ldots, n$. Hence, the smaller we make $g$ on $[0, X_{(n)})$, the smaller $\varphi_n(g)$ will become. It is not difficult to see that the smallest such $g$, with values of $g(X_i)$ fixed, must have the prescribed form. $\square$

Since $\widehat{h}_n$ is piecewise linear, it may be expressed as

$$\widehat{h}_n(t) = \widehat{a} + \sum_{j=1}^{k} \widehat{\nu}_j (\tau_j - t)_+ + \sum_{j=1}^{m} \widehat{\mu}_j (t - \eta_j)_+, \tag{2.1}$$

where $\widehat{\nu}_j, \widehat{a}, \widehat{\mu}_j \geq 0$. We let $\tau_j$ denote the points of change of slope of $\widehat{h}_n$ where $\widehat{h}_n$ is decreasing and let $\eta_j > 0$ denote the points of change of slope where $\widehat{h}_n$ is increasing. For simplicity, we assume that these are ordered. Also, we have $\tau_k \leq \eta_1$. As seen in the next lemma, the $\tau_j$'s and $\eta_j$'s correspond to "points of touch" or equality of processes defined on the one hand in terms of $\widehat{h}_n$ and the data, and on the other hand just in terms of the data. We therefore also refer to them as "touchpoints" repeatedly in the remainder of the paper.

It is convenient to define the MLE in terms of the minimization of the criterion function

$$\varphi_n(h) = \int_0^\infty \{H(t) - \log h(t) 1_{t \neq X_{(n)}}\} \, \mathrm{d}\mathbb{F}_n(t),$$

where $\mathbb{F}_n$ denotes the empirical distribution function of the data,

$$\mathbb{F}_n(t) = \frac{1}{n}\sum_{i=1}^{n} 1_{[0,t]}(X_i).$$

We will also use the notation $\mathbb{S}_n(t) = 1 - \mathbb{F}_n(t)$ for the empirical survival function.

**Lemma 2.2.** *Let $\widetilde{\mathbb{F}}_n(t) = (1/n)\sum_{i=1}^{n-1} 1_{[0,t]}(X_{(i)})$. A function $\widehat{h}_n$ minimizes $\varphi_n$ over $\mathcal{K}_+$ (and hence is the MLE) if and only if:*

$$\int_0^x \frac{x-t}{\widehat{h}_n(t)}\,\mathrm{d}\widetilde{\mathbb{F}}_n(t) = \int_0^x \int_0^t \mathbb{S}_n(s)\,\mathrm{d}s\,\mathrm{d}t \tag{2.2}$$

*for all $x \geq 0$ with equality at $\tau_i$ for $i = 1, \ldots, k$;*

$$\int_x^\infty \frac{t-x}{\widehat{h}_n(t)}\,\mathrm{d}\widetilde{\mathbb{F}}_n(t) = \int_x^\infty \int_t^\infty \mathbb{S}_n(s)\,\mathrm{d}s\,\mathrm{d}t \tag{2.3}$$

*for all $x \geq 0$ with equality at $\eta_j$ for $j = 1, \ldots, m$;*

$$\int_0^\infty \frac{1}{\widehat{h}_n(t)}\,\mathrm{d}\widetilde{\mathbb{F}}_n(t) \leq \int_0^\infty \mathbb{S}_n(t)\,\mathrm{d}t, \tag{2.4}$$

$$\int_0^\infty \widehat{H}_n(t)\,\mathrm{d}\mathbb{F}_n(t) = 1 - 1/n. \tag{2.5}$$

*Moreover, the minimizer $\widehat{h}_n$ satisfies*

$$\int_0^x \widehat{h}_n(t)\mathbb{S}_n(t)\,\mathrm{d}t = \mathbb{F}_n(x) \tag{2.6}$$

*for $x \in \{\tau_1, \ldots, \tau_k, \eta_1, \ldots, \eta_m\}$.*

**Remark 2.3.** As we assume a priori that $\widehat{h}_n(X_{(n)}) = \infty$, we may rewrite the left-hand side terms in (2.2)–(2.4) via

$$\int_A \frac{x-t}{\widehat{h}_n(t)}\,\mathrm{d}\widetilde{\mathbb{F}}_n(t) = \int_A \frac{x-t}{\widehat{h}_n(t)}1_{t \neq X_{(n)}}\,\mathrm{d}\mathbb{F}_n(t) = \int_A \frac{x-t}{\widehat{h}_n(t)}\,\mathrm{d}\mathbb{F}_n(t),$$

$$\int_A \frac{1}{\widehat{h}_n(t)}\,\mathrm{d}\widetilde{\mathbb{F}}_n(t) = \int_A \frac{1}{\widehat{h}_n(t)}1_{t \neq X_{(n)}}\,\mathrm{d}\mathbb{F}_n(t) = \int_A \frac{1}{\widehat{h}_n(t)}\,\mathrm{d}\mathbb{F}_n(t)$$

for any set $A$. We will hereafter use this latter formulation.

**Corollary 2.4.** *Let $\{\tau_i\}_{i=1}^k$ and $\{\eta_j\}_{j=1}^m$ denote the change points of $\widehat{h}_n$ as in (2.1). It follows that*

$$\int_0^{\tau_i} \frac{1}{\widehat{h}_n(t)}\,\mathrm{d}\mathbb{F}_n(t) = \int_0^{\tau_i} \mathbb{S}_n(u)\,\mathrm{d}u, \tag{2.7}$$

$$\int_{\eta_j}^{\infty} \frac{1}{\widehat{h}_n(t)} \, \mathrm{d}\mathbb{F}_n(t) = \int_{\eta_j}^{\infty} \mathbb{S}_n(u) \, \mathrm{d}u \tag{2.8}$$

for $i = 1, \ldots, k$ and $j = 1, \ldots, m$.

**Proof.** The function

$$\phi(x) \equiv \int_0^x \frac{x-t}{\widehat{h}_n(t)} \, \mathrm{d}\mathbb{F}_n(t) - \int_0^x \int_0^t \mathbb{S}_n(s) \, \mathrm{d}s \, \mathrm{d}t$$

is maximized at $\tau_i$ for $i = 1, \ldots, k$. Since it is also differentiable, (2.7) follows. A similar argument proves (2.8). □

**Proof of Lemma 2.2.** Consider any non-negative convex function $h$. It follows that there exists a non-negative constant $a$ and non-negative measures $\nu$ and $\mu$ (these measures have supports with intersection containing at most one point) such that

$$h(t) = a + \int_0^{\infty} (x-t)_+ \, \mathrm{d}\nu(x) + \int_0^{\infty} (t-x)_+ \, \mathrm{d}\mu(x).$$

For any function $\widehat{h}$ in $\mathcal{K}_+$, we calculate

$$\varphi_n(h) - \varphi_n(\widehat{h}) \geq \int_0^{\infty} \left\{ H(t) - \widehat{H}(t) + \left(1 - \frac{h(t)}{\widehat{h}(t)}\right) 1_{t \neq X_{(n)}} \right\} \mathrm{d}\mathbb{F}_n(t)$$

since $-\log x \geq 1 - x$. Plugging in the explicit form of $h$ from above, we find that the right-hand side is equal to

$$a \left\{ \int_{[0,\infty)} \left( t - \frac{1}{\widehat{h}(t)} 1_{t \neq X_{(n)}} \right) \mathrm{d}\mathbb{F}_n(t) \right\} + \left\{ \frac{n-1}{n} - \int_0^{\infty} \widehat{H}(t) \, \mathrm{d}\mathbb{F}_n(t) \right\}$$

$$+ \int_0^{\infty} \left\{ \int_0^x \int_0^t \mathbb{S}_n(s) \, \mathrm{d}s \, \mathrm{d}t - \int_0^x \frac{x-t}{\widehat{h}(t)} 1_{t \neq X_{(n)}} \, \mathrm{d}\mathbb{F}_n(t) \right\} \mathrm{d}\nu(x)$$

$$+ \int_0^{\infty} \left\{ \int_x^{\infty} \int_t^{\infty} \mathbb{S}_n(s) \, \mathrm{d}s \, \mathrm{d}t - \int_0^x \frac{t-x}{\widehat{h}(t)} 1_{t \neq X_{(n)}} \, \mathrm{d}\mathbb{F}_n(t) \right\} \mathrm{d}\mu(x).$$

This is non-negative if $\widehat{h}$ is a function which satisfies conditions (2.2)–(2.5). It follows that these conditions are sufficient to describe a minimizer of $\varphi_n$.

We next show that these conditions are necessary. To do this, we first define the directional derivative

$$\partial_\gamma \varphi_n(h) \equiv \lim_{\varepsilon \to 0} \frac{\varphi_n(h + \varepsilon \gamma) - \varphi_n(h)}{\varepsilon} = \int_0^{\infty} \left\{ \Gamma(t) - \frac{\gamma(t)}{h(t)} 1_{t \neq X_{(n)}} \right\} \mathrm{d}\mathbb{F}_n(t). \tag{2.9}$$

If $\widehat{h}_n$ minimizes $\varphi_n$, then for any $\gamma$ such that $\widehat{h}_n + \varepsilon\gamma$ is in $\mathcal{K}_+$ for sufficiently small $\varepsilon$, we must have $\partial_\gamma \varphi_n(\widehat{h}_n) \geq 0$. If, however, $\widehat{h}_n \pm \varepsilon\gamma$ is in $\mathcal{K}_+$ for sufficiently small $\varepsilon$, then $\partial_\gamma \varphi_n(\widehat{h}_n) = 0$.

If we choose, respectively, $\gamma(t) \equiv 1, (t-y)_+, (y-t)_+$, then $\widehat{h}_n + \varepsilon\gamma$ is in $\mathcal{K}_+$ and we obtain the inequalities in conditions (2.2)–(2.4). Since $(1 \pm \varepsilon)\widehat{h}_n$ is also in $\mathcal{K}_+$, for sufficiently small $\varepsilon$, we obtain (2.5). Choosing $\gamma = (\tau_i - t)_+, (t - \eta_j)_+$ yields the equalities in (2.2) and (2.3), respectively, since for each of these functions, $\widehat{h}_n \pm \varepsilon\gamma$ is in $\mathcal{K}_+$.

Lastly, we prove (2.6). For any $\tau_i$, define
$$\gamma(t) = \begin{cases} \widehat{h}_n(t) - \widehat{h}_n(\tau_i), & \text{for } t \in [0, \tau_i], \\ 0, & \text{otherwise.} \end{cases}$$

Since $(1 \pm \varepsilon)\gamma$ is also in $\mathcal{K}_+$, it follows that $\partial_\gamma \varphi_n(\widehat{h}_n) = 0$ and hence
$$0 = \left\{ \int_0^{\tau_i} \widehat{H}_n(t) \, d\mathbb{F}_n(t) - \mathbb{F}_n(\tau_i) + \widehat{H}_n(\tau_i) \mathbb{S}_n(\tau_i) \right\}$$
$$+ \widehat{h}_n(\tau_i) \left\{ \int_0^{\tau_i} \frac{1}{\widehat{h}_n(t)} \, d\mathbb{F}_n(t) - \int_0^{\tau_i} t \, d\mathbb{F}_n(t) - \tau_i \mathbb{S}_n(\tau_i) \right\}.$$

Integration by parts and Corollary 2.4 yield (2.6) for $x = \tau_i$. The case where $x = \eta_j$ is obtained in a similar manner, but using $\gamma(t) = (\widehat{h}_n(t) - \widehat{h}_n(\eta_j)) \mathbf{1}_{(\eta_j, \infty)}(t)$ and (2.5). □

The next corollary allows us to extend the equalities of the characterization of the MLE to some extra touchpoints. The significance of these equations will become clear in Section 5, where we consider asymptotics of the estimator.

**Corollary 2.5.** *Suppose that $\widehat{h}_n$ is strictly positive and recall the formulation given in (2.1). Then we also have that*

$$\int_0^{\eta_1} \frac{\eta_1 - t}{\widehat{h}_n(t)} \, d\mathbb{F}_n(t) = \int_0^{\eta_1} \int_0^s \mathbb{S}_n(u) \, du \, ds, \tag{2.10}$$

$$\int_{\tau_k}^\infty \frac{t - \tau_k}{\widehat{h}_n(t)} \, d\mathbb{F}_n(t) = \int_{\tau_k}^\infty \int_t^\infty \mathbb{S}_n(s) \, ds \, dt, \tag{2.11}$$

$$\int_0^{\tau_i} \frac{1}{\widehat{h}_n(t)} \, d\mathbb{F}_n(t) = \int_0^{\tau_i} \mathbb{S}_n(u) \, du, \tag{2.12}$$

$$\int_{\eta_j}^\infty \frac{1}{\widehat{h}_n(t)} \, d\mathbb{F}_n(t) = \int_{\eta_j}^\infty \mathbb{S}_n(u) \, du. \tag{2.13}$$

**Proof.** The first two equalities follow by noting that if $\widehat{h}_n$ is strictly positive, then for $\varepsilon$ sufficiently small, $\widehat{h}_n \pm \varepsilon\gamma$ is in $\mathcal{K}_+$ for $\gamma(t) = (t - \tau_k)_+, (\eta_1 - t)_+$. Arguing as for Corollary 2.4 proves the remaining identities. □

**Proposition 2.6.** *There exists a unique minimizer $\widehat{h}_n$ of $\varphi_n$ over $\mathcal{K}_+$.*

**Proof.** We will show that a minimizer exists by reducing the search to bounded positive convex functions on a compact domain. As this is a compact set, under the topology of uniform convergence, a minimizer of $\varphi_n$ exists (see [27], Theorems 10.6, 10.8 and 27.3).

We must first handle the issue of a compact domain. As we assume a priori that $\widehat{h}_n(X_{(n)}) = \infty$, we are really looking for the minimizer of the modified negative of the log-likelihood with domain $[0, X_{(n)})$. However, we have also argued that the minimizer must have the specific functional form as described in Proposition 2.1. Therefore, it is sufficient to reduce the domain to $[0, X_{(n-1)} + \delta]$, for any $\delta > 0$, since $\widehat{h}_n$ is then extended linearly beyond $X_{(n-1)} + \delta$ in a unique manner. It will therefore be sufficient to show that we may reduce the search to functions bounded on $[0, X_{(n-1)}]$, with a derivative at $X_{(n-1)}$ which is bounded above.

Recall that the minimizer must satisfy (2.5). We therefore reduce our search to the class of functions which satisfy this condition. For any such $h$, write $h = h_+ + h_-$, where $h_+$ is increasing and $h_-$ is decreasing. It follows that for any $x$,

$$1 \geq \int_0^\infty H(t) \, \mathrm{d}\mathbb{F}_n(t) = \int_0^\infty h(t) \mathbb{S}_n(t) \, \mathrm{d}t \geq h_-(x) \int_0^x \mathbb{S}_n(t) \, \mathrm{d}t.$$

A similar bound for $h_+$ yields

$$h(x) \leq \left( \int_0^x \mathbb{S}_n(t) \, \mathrm{d}t \right)^{-1} + \left( \int_x^\infty \mathbb{S}_n(t) \, \mathrm{d}t \right)^{-1} \equiv M_n(x) \qquad (2.14)$$

for all $x$ in $(0, X_{(n)})$. Thus we know that $h(x)$ must be bounded for $x \in (0, X_{(n-1)}]$.

To show that $h$ is also bounded at zero, we need to show that $h'(X_{(1)})$ is bounded from below. Assuming that it is negative, we may write, for $0 < x \leq X_{(1)}$,

$$h(X_{(1)}) + h'(X_{(1)})(x - X_{(1)}) = h(x) \leq M_n(x).$$

Fixing $x^* > 0$ and less than $X_{(1)}$, we then obtain that

$$h'(X_{(1)}) \geq (M_n(x^*) - h(X_{(1)}))/(x^* - X_{(1)}),$$

from which it follows that $h$ must be bounded on the set $[0, X_{(n-1)}]$.

By (2.5), we also have that

$$n \geq H(X_{(n)}) \geq \int_{X_{(n-1)}}^{X_{(n-1)}+\delta} h(t) \, \mathrm{d}t = \int_{X_{(n-1)}}^{X_{(n-1)}+\delta} \{h(X_{(n-1)}) + h'(X_{(n-1)})(t - X_{(n-1)})\} \, \mathrm{d}t$$

if $h$ is increasing on $[X_{(n-1)}, X_{(n)})$. This implies that $h'(X_{(n-1)})$ is bounded above, completing the proof.

We now show uniqueness. Suppose that $h_1$ and $h_2$ both minimize $\varphi_n$. Then, by (2.5), $\varphi_n(h_1)$ and $\varphi_n(h_2)$ differ only in the term $-\int_0^\infty \log h_i(t) \mathbf{1}(t \neq X_{(n)}) \, \mathrm{d}\mathbb{F}_n(t)$. However, this term is strictly convex and it follows that $h_1(X_{(i)}) = h_2(X_{(i)})$ for all $i = 1, \ldots, n-1$.

Let $\bar{h} = (h_1 + h_2)/2$. By linearity, we have that $\varphi_n(h_1) = \varphi(h_2) = \varphi(\bar{h})$, which implies that $\bar{h}$ is also a minimizer. However, the only way that this is possible is if $\bar{h}$ also satisfies the conditions of Proposition 2.1. This implies that one of the following holds:

1. Both $h_1$ and $h_2$ are increasing and $h_1(0) = h_2(0) = 0$. In this case, they must have the same locations for their changes of slope, as otherwise $\bar{h}$ violates Proposition 2.1.
2. Point 1 above does not hold. Then, by the same argument as above, if $h_1$ and $h_2$ have at least one change of slope in an interval between observations (or between zero and $X_{(1)}$), then these locations of change of slope must be equal.

If the first case holds, then it is not difficult to see that $h_1 \equiv h_2$ on $[0, X_{(n)})$, as $h_1(t) = h_1(t) = 0$ on $[0, \tau_1]$ and $h_1(X_i) = h_2(X_i)$ for all observation points.

In the second case, we use a different argument. We know that neither $h_1$ nor $h_2$ have touchpoints before $X_{(1)}$. Let $t^*$ denote the first touchpoint of $h_1$ and (without loss of generality) assume that the first touchpoint of $h_2$ is greater than $t^*$. Hence, by (2.6),

$$h_1(X_{(1)}) = h_2(X_{(1)}), \qquad \int_0^{t^*} h_1(t)\,dt = \mathbb{F}_n(t^*).$$

Now, $\bar{h} = (h_1 + h_2)/2$ and $h_2$ are also minimizers of the MLE criterion function $\varphi_n$. Also, $\bar{h}$ has a touchpoint at $t^*$ and $\bar{h}(X_{(1)}) = h_2(X_{(1)})$.

Averaging $\bar{h}$ with $h_2$ yields the functions $\bar{h}_l = 2^{-l}(h_1 - h_2) + h_2$, which satisfy

$$\bar{h}_l(X_{(1)}) = h_2(X_{(1)}), \qquad \int_0^{t^*} \bar{h}_l(t)\,dt = \mathbb{F}_n(t^*) \qquad \text{for all } l \geq 1.$$

Since $\bar{h}_l \to h_2$ pointwise, it follows from the dominated convergence theorem that $\int_0^{t^*} h_2(t)\,dt = \mathbb{F}_n(t^*)$. Therefore, since $h_1$ and $h_2$ are both linear on $[0, t^*]$ with $h_1(X_{(1)}) = h_2(X_{(1)})$ and $\int_0^{t^*} h_1(t)\,dt = \int_0^{t^*} h_2(t)\,dt$, they must have both the same value *and slope* at $X_{(1)}$. That is, both $h_1(X_{(1)}) = h_2(X_{(2)})$ and $h_1'(X_{(1)}) = h_2'(X_{(2)})$ hold.

Now, write

$$h_1(t) = a_1 + b_1 t + \sum_{i=1}^{m_1 - 1} \nu_{i,1}(t - t_{i,1})_+,$$

$$h_2(t) = a_2 + b_2 t + \sum_{i=1}^{m_2 - 1} \nu_{i,2}(t - t_{i,2})_+,$$

where $X_{(1)} < t_{1,j} < t_{2,j} < \cdots < t_{m_j - 1, j} < X_{(n)}$, $j = 1, 2$, and where $h_1(X_{(i)}) = h_2(X_{(i)})$ for $i = 1, \ldots, n$. We also assume that $\nu_{i,j} > 0$ for $i = 1, \ldots, m_j - 1$, $j = 1, 2$. This implies, in particular, that $h_j(t) = a_j + b_j t$ for $t \leq t_{1,j}$, $j = 1, 2$, and since $X_{(1)} < t_{1,j}$, $j = 1, 2$, $h_1(X_{(1)}) = h_2(X_{(1)})$. Thus $a_1 + b_1 X_{(1)} = a_2 + b_2 X_{(1)}$. From the argument above, it follows that $b_1 = h_1'(X_{(1)}) = h_2'(X_{(1)}) = b_2$. We conclude that $a_1 = a_2$ and $b_1 = b_2$ so that $h_1(t) = h_2(t)$ for $0 \leq t \leq t^*$. It also follows that $t_{1,1} = t_{1,2}$.

Repeating this argument on the interval $[t^*, t^{**}]$ with $t^{**} = \min\{t_{2,1}, t_{2,2}\}$ shows that $\nu_{1,1} = \nu_{1,2}$ or $t_{2,1} = t_{2,2}$. Proceeding by induction yields $\nu_{j,1} = \nu_{j,2}$ and $t_{j+1,1} = t_{j+1,2}$ for $j = 1, \ldots, m_1 - 1 = m_2 - 1$, hence uniqueness. □

## 3. Consistency

**Theorem 3.1.** *Suppose that $X_1, \ldots, X_n$ are i.i.d. random variables with convex hazard function $h_0$ and corresponding distribution function $F_0$. Let $T_0 \equiv T_0(F_0) \equiv \inf\{t : F_0(t) = 1\}$. The MLE $\widehat{h}_n(t)$ is then consistent for all $t \in (0, T_0)$. Also, for all $\delta > 0$,*

$$\sup_{\delta \leq t \leq T_0 - \delta} |\widehat{h}_n(t) - h(t)| \to 0, \qquad \textit{almost surely},$$

*if $T_0 < \infty$. If $T_0 = \infty$, the above statement holds with $T_0 - \delta$ replaced by any $K < \infty$.*

*Remark 3.2.* If $h_0$ is increasing at 0, then one can show that $\widehat{h}_n$ is not consistent at zero. This is a frequently occurring difficulty of shape constrained estimators; see, for example, [1, 12, 30].

**Proof.** We first show that $\widehat{h}_n$ is bounded appropriately so that we can select convergent subsequences. Decompose $\widehat{h}_n$ into its decreasing and increasing parts: $\widehat{h}_n = \widehat{h}_{n,\downarrow} + \widehat{h}_{n,\uparrow}$. Then, arguing as in (2.14), it follows from (2.5) that

$$\widehat{h}_{n,\downarrow}(x) \leq \frac{1}{\int_0^x \mathbb{S}_n(t)\,dt}, \tag{3.1}$$

where the right-hand side is almost surely bounded and, in fact, converges almost surely to $1/\int_0^x S_0(t)\,dt < \infty$ for all $x > 0$. Also,

$$\widehat{h}_{n,\uparrow}(x) \leq \frac{1}{\int_x^{x+\delta} \mathbb{S}_n(t)\,dt}, \tag{3.2}$$

where the right-hand side is almost surely bounded for $x \in (\mathrm{supp}(F_0))^\circ$ and converges almost surely to $1/\int_x^{x+\delta} S_0(t)\,dt < \infty$.

Now, take $\gamma = h_0$ in the directional derivative (2.9). It follows that

$$0 \leq \lim_{\varepsilon \downarrow 0} \frac{\varphi_n(\widehat{h}_n + \varepsilon h_0) - \varphi_n(\widehat{h}_n)}{\varepsilon} = \int_0^\infty \left\{ H_0(t) - \frac{h_0(t)}{\widehat{h}_n(t)} \right\} d\mathbb{F}_n(t),$$

noting that $\widehat{h}_n(X_{(n)}) = \infty$, and hence

$$\int_0^\infty \frac{h_0(t)}{\widehat{h}_n(t)} d\mathbb{F}_n(t) \leq \int_0^\infty H_0(t)\,d\mathbb{F}_n(t) \underset{\text{a.s.}}{\to} \int_0^\infty H_0(t)\,dF_0(t) = 1.$$

Fix any $0 < a < b < \infty$ such that $a, b \in (\text{supp}(F_0))^\circ$. It follows that $\lim_n X_{(n)} > b$ with probability one (this can be shown using the Borel–Cantelli theorem). Also, $\sup |\mathbb{F}_n(t) - F_0(t)| \to_{\text{a.s.}} 0$ by the Glivenko–Cantelli lemma. Both of these events occur on the set $\Omega$, with $P(\Omega) = 1$. Fix $\omega \in \Omega$. We will show that $\widehat{h}_n \to h_0$ for such an $\omega$.

Let $\{n'\}$ denote any subsequence of $\{n\}$. By the bounds in (3.1) and (3.2) (which are finite for our choice of $\omega$), using a classical diagonalization argument and the continuity of convex functions, we may extract a further subsequence $\{n''\}$ such that $\widehat{h}_{n''} \to \widehat{h}$ pointwise on $[a, b]$, where the limit $\widehat{h}$ must be convex. We denote this subsequence $\{n\}$ to simplify notation.

From Fatou's lemma, it follows that

$$\int_a^b \frac{h_0^2(t)}{\widehat{h}_n(t)} S_0(t)\,dt = \int_a^b \frac{h_0(t)}{\widehat{h}_n(t)} f_0(t)\,dt \leq \liminf_n \int_a^b \frac{h_0(t)}{\widehat{h}_n(t)}\,d\mathbb{F}_n(t)$$
$$\leq \limsup_n \int_0^\infty \frac{h_0(t)}{\widehat{h}_n(t)}\,d\mathbb{F}_n(t) \leq \lim_n \int_0^\infty H_0(t)\,d\mathbb{F}_n(t) \leq 1.$$

Note that this implies that if $\widehat{h}(t) = 0$, then $h_0(t) = 0$. By (2.5) and integration by parts, we see that $1 \geq \int_{[0, X_{(n)})} \widehat{h}_n(t) \mathbb{S}_n(t)\,dt$. Therefore, again applying Fatou's lemma,

$$1 \geq \int_a^b \widehat{h}(t) S_0(t)\,dt.$$

It also follows that

$$0 \leq \int_a^b \frac{(\widehat{h}(t) - h_0(t))^2}{\widehat{h}(t)} S_0(t)\,dt$$
$$= \int_a^b \widehat{h}(t) S_0(t)\,dt - 2 \int_a^b h_0(t) S_0(t)\,dt + \int_a^b \frac{h_0^2(t)}{\widehat{h}(t)} S_0(t)\,dt$$
$$\leq 2 - 2 \int_a^b h_0(t) S_0(t)\,dt.$$

Define $\widehat{h} = h_0$ for $t \notin [a, b]$, which allows us to let both $1/a$ and $b \to \infty$ in the above display. Since $\int_0^\infty h_0(t) S_0(t)\,dt = 1$, it follows that

$$\int_0^\infty \frac{(\widehat{h}(t) - h_0(t))^2}{\widehat{h}(t)} S_0(t)\,dt = 0$$

and this implies that $\widehat{h}(t) = h_0(t)$ for all $t \in [a, b]$.

We have thus shown that every subsequence $\{\widehat{h}_n(x)\}$ has a further subsequence which converges to the true hazard function $h_0(x)$ pointwise, for all $x \in (\text{supp}\, F_0)^\circ$. It follows that $\{\widehat{h}_n\}$ converges to $h_0$ pointwise. By Theorem 10.8, page 90, [27], this implies that the claimed uniform convergence on $[a,b]$ also holds. As this happens for any $\omega \in \Omega$, and $P(\Omega) = 1$, we have proven the result. □

**Corollary 3.3.** *Suppose that $h_0''$ is continuous and strictly positive at $x_0 \in (\operatorname{supp} F_0)^\circ$. It follows that there exist touchpoints $\tau_n \leq x_0 \leq \eta_n$ such that $\tau_n, \eta_n \to x_0$ in probability.*

**Proof.** Let $\eta_n, \tau_n$ be touchpoints such that $\tau_n \leq x_0 \leq \eta_n$. If $\tau_n$ does not exist, then set $\tau_n = 0$, and $\eta_n = \infty$ otherwise. Suppose that it is not the case that $\tau_n, \eta_n \to_p x_0$. It then follows from Theorem 3.1 that there exists an interval $I = [a,b]$ such that $x_0 \in I$ for $|I| > 0$, $\limsup_n \tau_n \leq a$ and $\liminf_n \eta_n \geq b$ almost surely, and, lastly, $\lim \widehat{h}_n(t) \to_{\text{a.s.}} h_0(t)$ on $I$. However, this implies that $h_0(t)$ is linear on $I$, which is a contradiction. □

From consistency of the estimator, we also obtain consistency of the derivatives.

**Corollary 3.4.** *Suppose that $x \in (a,b)$ and $\sup_{a \leq t \leq b} |\widehat{h}_n(t) - h_0(t)| \to_{\text{a.s.}} 0$. Then $\widehat{h}_n'(x) \to_{\text{a.s.}} h_0'(x)$ at all continuity points $x$ of $h_0'$ on $(a,b)$.*

This follows from the following simple result for convex functions, proved in [12].

**Lemma 3.5.** *Suppose that $\bar{h}_n$ is a sequence of convex functions satisfying $\sup_{a \leq x \leq b} |\bar{h}_n(t) - h_0(t)| \to 0$ with probability one. Then (also with probability one), for all $x \in (a,b)$,*

$$-\infty < h_0'(x^-) \leq \liminf_{n \to \infty} \bar{h}_n'(x^-) \leq \limsup_{n \to \infty} \bar{h}_n'(x^+) \leq h_0'(x^+) < \infty.$$

## 4. Asymptotic lower bounds for the minimax risk

Define the class of densities $\mathcal{C}$ by

$$\mathcal{C} = \bigg\{ f : [0, \infty) \to [0, \infty) \colon \int_0^\infty f(x) \, dx = 1,$$

$$h(x) = f(x)/(1 - F(x)) \text{ is convex}, h(x) > 0 \text{ for all } x > 0 \bigg\}.$$

We want to derive asymptotic lower bounds for the local minimax risks for estimating the convex hazard function $h$ and its derivative at a fixed point. The $L_1$-minimax risk for estimating a functional $T$ of $f_0$ based on a sample $X_1, \ldots, X_n$ of size $n$ from $f_0$ which is known to be in a subset $\mathcal{C}_n$ of $\mathcal{C}$ is defined by

$$MMR_1(n, T, \mathcal{C}_n) = \inf_{T_n} \sup_{f \in \mathcal{C}_n} E_f |T_n - Tf|, \tag{4.1}$$

where the infimum ranges over all possible measurable functions $T_n = t_n(X_1, \ldots, X_n)$ mapping $\mathbb{R}^n$ to $\mathbb{R}$. The shrinking classes $\mathcal{C}_n$ used here are Hellinger balls centered at $f_0$,

$$\mathcal{C}_{n,\tau} = \left\{ f \in \mathcal{C} : H^2(f, f_0) \equiv \frac{1}{2} \int_0^\infty (\sqrt{f(z)} - \sqrt{f_0(z)})^2 \, dz \leq \tau/n \right\}.$$

Consider estimation of

$$T_1(f) = \frac{f(x_0)}{1 - F(x_0)} = h(x_0), \qquad T_2(f) = h'(x_0). \tag{4.2}$$

Let $f_0 \in \mathcal{C}$ and $x_0 > 0$ be fixed such that $h_0$ is twice continuously differentiable at $x_0$. Define, for $\varepsilon > 0$, the functions $h_\varepsilon$ as follows:

$$h_\varepsilon(z) = \begin{cases} h_0(x_0 - \varepsilon c_\varepsilon) + (z - x_0 + \varepsilon c_\varepsilon) h_0'(x_0 - \varepsilon c_\varepsilon), & z \in [x_0 - \varepsilon c_\varepsilon, x_0 - \varepsilon], \\ h_0(x_0 + \varepsilon) + (z - x_0 - \varepsilon) h_0'(x_0 + \varepsilon), & z \in [x_0 - \varepsilon, x_0 + \varepsilon], \\ h_0(z), & \text{otherwise.} \end{cases}$$

Here, $c_\varepsilon$ is chosen so that $h_\varepsilon$ is continuous at $x_0 - \varepsilon$. Using continuity of $h_\varepsilon$ and a second order expansion of $h_0$, it follows that $c_\varepsilon = 3 + o(1)$ as $\varepsilon \to 0$. Now, define $f_\varepsilon$ by

$$f_\varepsilon(z) = \exp(-H_\varepsilon(z)) h_\varepsilon(z),$$

where $H_\varepsilon(z) \equiv \int_0^z h_\varepsilon(u) \, du$. It follows easily that

$$T_1(f_\varepsilon) - T_1(f_0) = \tfrac{1}{2} h_0''(x_0) \varepsilon^2 + o(\varepsilon^2), \tag{4.3}$$

$$T_2(f_\varepsilon) - T_2(f_0) = h_0''(x_0) \varepsilon + o(\varepsilon). \tag{4.4}$$

Furthermore, the following lemma holds.

**Lemma 4.1.** *Under the above assumptions,*

$$H^2(f_\varepsilon, f_0) = \frac{2}{5} \frac{h_0''(x_0)^2 (1 - F(x_0))}{h_0(x_0)} \varepsilon^5 + o(\varepsilon^5) \equiv \nu_0 \varepsilon^5 + o(\varepsilon^5).$$

**Proof.** The lemma follows from Lemma 2 of [18] and

$$\int \frac{(f_\varepsilon(x) - f_0(x))^2}{f_0(x)} \, dx = \frac{16}{5} \frac{h_0''(x_0)^2 (1 - F(x_0))}{h_0(x_0)} \varepsilon^5 + o(\varepsilon^5).$$

This is achieved by careful calculation. □

Combining (4.3) and (4.4) with the lemma, it follows that

$$|T_1(f_{(\varepsilon/\nu_0)^{1/5}}) - T_1(f_0)| \geq \left( \frac{h_0(x_0) \sqrt{h_0''(x_0)}}{S_0(x_0) 8 \sqrt{2}} \right)^{2/5} \varepsilon^{2/5} (1 + o(1)),$$

$$|T_2(f_{(\varepsilon/\nu_0)^{1/5}}) - T_2(f_0)| \geq \left( \frac{5 h_0(x_0) h_0''(x_0)^3}{2 S_0(x_0)} \right)^{1/5} \varepsilon^{1/5} (1 + o(1)).$$

From these calculations, together with Lemma 5.1 of [12], we have the following result. Along with Theorem 1.2, it indicates that $\widehat{h}_n(x_0)$ and $\widehat{h}'_n(x_0)$ achieve optimal rates and also have the correct dependence on the parameters $h''(x_0)$ and $h(x_0)$ (up to absolute constants).

**Theorem 4.2 (Minimax risk lower bound).** *For the functionals $T_1$ and $T_2$ as defined in (4.2), and with $MMR_1(n, T, \mathcal{C}_{n,\tau})$ as defined in (4.1),*

$$\sup_{\tau>0} \limsup_{n\to\infty} n^{2/5} MMR_1(n, T_1, \mathcal{C}_{n,\tau}) \geq \frac{1}{4}\left(\frac{h_0(x_0)\sqrt{h_0''(x_0)}}{S_0(x_0)\mathrm{e}8\sqrt{2}}\right)^{2/5} \quad and$$

$$\sup_{\tau>0} \limsup_{n\to\infty} n^{1/5} MMR_1(n, T_2, \mathcal{C}_{n,\tau}) \geq \frac{1}{4}\left(\frac{1}{4\mathrm{e}}\frac{h_0(x_0)h_0''(x_0)^3}{2S_0(x_0)}\right)^{1/5}.$$

In particular, Theorem 1.2 shows that the MLE achieves the optimal pointwise rate of convergence, $n^{2/5}$, at points $x_0$ with $h''(x_0) > 0$. Convergence rates over the larger class of bathtub-shaped functions would be slower: the MLE of a $U$-shaped hazard is known to converge locally at rate $n^{1/3}$; see, for example, [2].

## 5. Rates of convergence

In this section, we identify the local rates of convergence of the MLE. Fix a point $x_0 \in (\operatorname{supp} F_0)^\circ$. To obtain the results, we assume that $h_0''(\cdot)$ is continuous and strictly positive in a neighborhood of $x_0$ and that $h(x_0) > 0$.

### 5.1. Some useful estimates

For $0 < x \leq y$, define

$$U_n(x, y) = \int_x^y \left\{\frac{z - (1/2)(x+y)}{\widehat{h}_n(z)}\right\} \mathrm{d}(\mathbb{F}_n - F_0)(z).$$

**Lemma 5.1.** *Let $x_0 \in (\operatorname{supp} F_0)^\circ$. Then, for each $\varepsilon > 0$, there exist constants $\delta, c_0, n_0$ and (positive) random variables $M_n$ (independent of $x, y$), of order $\mathrm{O}_p(1)$, such that for each $|x - x_0| < \delta$,*

$$|U_n(x, y)| \leq \varepsilon(y-x)^4 + n^{-4/5} M_n, \qquad 0 \leq y - x \leq c_0, \tag{5.1}$$

*for all $n \geq n_0$.*

**Proof.** Note that $U_n = (\mathbb{P}_n - P_0)(g_{x,y,\widehat{h}_n})$, where

$$g_{x,y,h}(z) \equiv \frac{f_{x,y}(z)}{h(z)} 1_{[x,y]}(z)$$

and, in view of the consistency established in Theorem 3.1, $\widehat{h}_n$ is a convex function uniformly close to $h_0$ on neighborhoods of $x_0$. This leads to the consideration of the class of functions

$$\mathcal{F}_{x,R} \equiv \left\{ \begin{array}{c} z \mapsto g_{x,y,h}(z) : x \leq y \leq x+R, h \text{ convex,} \\ \|h - h_0\|_{x_0-\delta}^{x_0+\delta+c_0} \leq \gamma \end{array} \right\}$$

with $\gamma \equiv \inf_{x_0-\delta \leq x \leq x_0+\delta+c_0} h_0(x)/2$ and we define $G_n \equiv \{\|\widehat{h}_n - h_0\|_{x_0-\delta}^{x_0+\delta+c_0} \leq \gamma\}$. The class $\mathcal{F}_{x,R}$ has an envelope function $F_{x,R}(z) = \gamma^{-1}\{(z-x)1_{[x,x+R]}(z) + 2^{-1}R1_{[x,x+R]}(z)\}$ and hence the following second moment bound holds:

$$E\{[F_{x,R}]^2\} = \frac{1}{\gamma^2}\int_{[x,x+R]}[(z-x)+R/2]^2 f_0(z)\,dz \leq \frac{13}{12\gamma^2}\|f_0\|_{x_0-\delta}^{x_0+\delta}R^3.$$

Furthermore, $\log N_{[]}(\varepsilon, \mathcal{F}_{x,R}, L_2(P_0)) \leq K/\varepsilon^{1/2}$ for some constant $K$ by [29], Theorem 2.7.10, page 159, and a straightforward bracketing argument. It then follows from [29], Theorems 2.14.2 and 2.14.5, pages 240 and 244, that

$$E\left\{\left(\sup_{f \in \mathcal{F}_{x,R}} |(\mathbb{P}_n - P_0)(f)|\right)^2\right\} \leq \frac{1}{n}K'E\{[F_{x,R}(X_1)]^2\} = O(n^{-1}R^3). \qquad (5.2)$$

Define $M_n(\omega)$ as the infimum (possibly $+\infty$) of those values such that (5.1) holds and define $A(n,j)$ to be the set $[(j-1)n^{-1/5}, jn^{-1/5})$. Then, for $m$ constant,

$$\begin{aligned} P(M_n > m) &\leq P([M_n > m] \cap G_n) + P(G_n^c) \\ &\leq P([\exists u : |U_n(x, x+u)| > \varepsilon u^4 + n^{-4/5}m] \cap G_n) + P(G_n^c) \\ &\leq \sum_{j \geq 1} P([\exists u \in A(n,j) : n^{4/5}|U_n(x, x+u)| > \varepsilon(j-1)^4 + m] \cap G_n) + P(G_n^c). \end{aligned}$$

The $j$th summand is hence bounded by

$$n^{8/5}E\left[\sup_{u \in A(n,j)} |U_n(x, x+u)m|^2 1_{G_n}\right]/[m+\varepsilon(j-1)^4]^2 \leq C\frac{j^3}{[m+\varepsilon(j-1)^4]^2}$$

due to (5.2). Thus it follows, using Theorem 3.1 to conclude that $P(G_n^c) \to 0$, that

$$\limsup_{n \to \infty} P(M_n > m) \leq C \sum_{j=1}^{\infty} \frac{j^3}{[m+\varepsilon(j-1)^4]^2},$$

where the sum in the bound is finite and converges to zero as $m \to \infty$. This completes the proof of the claim. □

A similar approach proves the following for the function

$$V_n(x,y) = \int_x^y \{z - (x+y)/2\}(\mathbb{S}_n(z) - S_0(z))\,dz.$$

**Lemma 5.2.** *Let $x_0 \in (\operatorname{supp} F_0)^\circ$. Then, for each $\varepsilon > 0$, there exist constants $\delta, c_0 > 0$ and (positive) random variables $M_n$ (independent of $x, y$) of order $\mathrm{O}_p(1)$ such that for each $|x - x_0| < \delta$,*

$$|V_n(x,y)| \leq \varepsilon n^{-1/5}(y-x)^4 + n^{-1} M_n, \qquad 0 \leq y - x \leq c_0. \tag{5.3}$$

### 5.2. Asymptotic behavior of touchpoints and resulting bounds

**Lemma 5.3.** *Let $x_0 > 0$ be a point at which $h_0$ has a continuous and strictly positive second derivative, and where $h(x_0) > 0$. Let $\xi_n$ be any sequence of numbers converging to $x_0$ and define $\tau_n$ and $\eta_n$ to be the largest touchpoint of $\widehat{h}_n$ smaller than $\xi_n$ and the smallest touchpoint larger than $\xi_n$, respectively. Then*

$$\eta_n - \tau_n = \mathrm{O}_p(n^{-1/5}).$$

**Proof.** By Theorem 3.1, we know that $\widehat{h}_n$ is positive near $x_0$ for large enough $n$. Also, it is either strictly increasing or strictly decreasing in a neighborhood of $x_0$, or it is locally flat. If $\widehat{h}_n$ is decreasing between $\tau_n$ and $\eta_n$, then (2.7) and (2.2) with equality at both $\eta_n$ and $\tau_n$ hold. If, instead, $\widehat{h}_n$ is increasing, then (2.8) and (2.3) with equality at both $\eta_n$ and $\tau_n$ hold. There is only the potential for a problem in the locally flat case. However, since $\widehat{h}_n$ is strictly positive, by Corollary 2.5, we can extend the necessary equalities to this case as well. Therefore, we need only consider two cases, $\widehat{h}_n$ is either non-increasing or non-decreasing on $[\tau_n, \eta_n]$.

We first assume that $\widehat{h}_n$ is non-increasing on $[\tau_n, \eta_n]$. Define

$$\widehat{\mathcal{H}}_{n,\downarrow}(z) = \int_0^z \frac{z-t}{\widehat{h}_n(t)}\,\mathrm{d}\mathbb{F}_n(t) \quad \text{and} \quad \mathbb{A}_{n,\downarrow}(z) = \int_0^z \mathbb{S}_n(t)\,\mathrm{d}t, \tag{5.4}$$

and let $m_n$ be the midpoint of $[\tau_n, \eta_n]$, $m_n = (\tau_n + \eta_n)/2$. We may then calculate

$$\widehat{\mathcal{H}}_{n,\downarrow}(m_n) = \int_{m_n}^{\eta_n} \frac{x - m_n}{\widehat{h}_n(x)}\,\mathrm{d}\mathbb{F}_n(x) + \widehat{\mathcal{H}}_{n,\downarrow}(\eta_n) - (\eta_n - m_n)\widehat{\mathcal{H}}'_{n,\downarrow}(\eta_n)$$

$$= \int_{\tau_n}^{m_n} \frac{m_n - x}{\widehat{h}_n(x)}\,\mathrm{d}\mathbb{F}_n(x) + \widehat{\mathcal{H}}_{n,\downarrow}(\tau_n) + (m_n - \tau_n)\widehat{\mathcal{H}}'_{n,\downarrow}(\tau_n).$$

From (2.2), we know that $2\widehat{\mathcal{H}}_{n,\downarrow}(m_n) \leq 2\int_0^{m_n} \mathbb{A}_{n,\downarrow}(t)\,\mathrm{d}t$. The equality in (2.2), together with (2.7), allows us to rewrite this as $0 \geq L_{1,\downarrow} + L_{2,\downarrow}$, where $L_{1,\downarrow}$ is equal to

$$\int_{m_n}^{\eta_n} \frac{x - m_n}{\widehat{h}_n(x)}\,\mathrm{d}\mathbb{F}_n(x) + \int_{\tau_n}^{m_n} \frac{m_n - x}{\widehat{h}_n(x)}\,\mathrm{d}\mathbb{F}_n(x) - \frac{\eta_n - \tau_n}{4}\{\widehat{\mathcal{H}}'_{n,\downarrow}(\eta_n) - \widehat{\mathcal{H}}'_{n,\downarrow}(\tau_n)\}$$

$$= \int_{m_n}^{\eta_n} \frac{x - (1/2)(\eta_n + m_n)}{\widehat{h}_n(x)}\,\mathrm{d}\mathbb{F}_n(x) + \int_{\tau_n}^{m_n} \frac{(1/2)(\tau_n + m_n) - x}{\widehat{h}_n(x)}\,\mathrm{d}\mathbb{F}_n(x)$$

and

$$L_{2,\downarrow} = \int_{m_n}^{\eta_n} \mathbb{A}_{n,\downarrow}(x)\,\mathrm{d}x - \int_{\tau_n}^{m_n} \mathbb{A}_{n,\downarrow}(x)\,\mathrm{d}x - \frac{1}{4}(\eta_n - \tau_n)\{\mathbb{A}_{n,\downarrow}(\eta_n) - \mathbb{A}_{n,\downarrow}(\tau_n)\}$$

$$= -\left\{\int_{m_n}^{\eta_n}\left\{x - \frac{1}{2}(\eta_n + m_n)\right\}\mathbb{S}_n(x)\,\mathrm{d}x + \int_{\tau_n}^{m_n}\left\{\frac{1}{2}(\tau_n + m_n) - x\right\}\mathbb{S}_n(x)\,\mathrm{d}x\right\},$$

by integration by parts.

Now replace $\mathbb{F}_n$ by the true $F_0$ in the definition of $L_{1,\downarrow}$ to obtain

$$L_{1,\downarrow}^0 \equiv \int_{m_n}^{\eta_n} \frac{x - (1/2)(\eta_n + m_n)}{\widehat{h}_n(x)}\,\mathrm{d}F_0(x) + \int_{\tau_n}^{m_n} \frac{(1/2)(\tau_n + m_n) - x}{\widehat{h}_n(x)}\,\mathrm{d}F_0(x)$$

$$= \int_{m_n}^{\eta_n}\left\{x - \frac{1}{2}(\eta_n + m_n)\right\}\left\{\frac{1}{\widehat{h}_n(x)} - \frac{1}{h_0(x)}\right\}\mathrm{d}F_0(x)$$

$$+ \int_{\tau_n}^{m_n}\left\{\frac{1}{2}(\tau_n + m_n) - x\right\}\left\{\frac{1}{\widehat{h}_n(x)} - \frac{1}{h_0(x)}\right\}\mathrm{d}F_0(x) - L_{2,\downarrow}^0,$$

where

$$L_{2,\downarrow}^0 = -\int_{m_n}^{\eta_n}\left\{x - \frac{1}{2}(\eta_n + m_n)\right\}S_0(x)\,\mathrm{d}x - \int_{\tau_n}^{m_n}\left\{\frac{1}{2}(\tau_n + m_n) - x\right\}S_0(x)\,\mathrm{d}x.$$

Next, using a Taylor expansion of order 2 on the function $1/\widehat{h}_n(x) - 1/h_0(x)$ and about the point $m_n$, we obtain

$$L_{1,\downarrow}^0 + L_{2,\downarrow}^0 = \frac{1}{192}\left\{\frac{h_0''(x_0)}{h_0^2(x_0)}f_0(x_0)\right\}(\eta_n - \tau_n)^4 + \mathrm{o}((\eta_n - \tau_n)^4) \tag{5.5}$$

since both $\widehat{h}_n$ and $\widehat{h}_n'$ are consistent by Theorem 3.1, $\widehat{h}_n''(x) = 0$ on $(\tau_n, \eta_n)$ and because $\tau_n - \eta_n = \mathrm{o}_p(1)$ by Corollary 3.3. Therefore, by Lemmas 5.1 and 5.2, together with the above calculations, we may write

$$0 \geq L_{1,\downarrow} + L_{2,\downarrow}$$
$$= L_{1,\downarrow}^0 + L_{2,\downarrow}^0 + (L_{1,\downarrow} - L_{1,\downarrow}^0) + (L_{2,\downarrow} - L_{2,\downarrow}^0)$$
$$\geq L_{1,\downarrow}^0 + L_{2,\downarrow}^0 - \varepsilon(\eta_n - \tau_n)^4 - \mathrm{O}_p(n^{-4/5}) - \varepsilon n^{-1/5}(\eta_n - \tau_n)^4 - \mathrm{O}_p(n^{-1})$$
$$= \frac{1}{192}\left\{\frac{h_0''(x_0)}{h_0^2(x_0)}f(x_0) - 192\varepsilon\right\}(\eta_n - \tau_n)^4 + \mathrm{o}((\eta_n - \tau_n)^4) - \mathrm{O}_p(n^{-4/5}).$$

We choose $\varepsilon$ sufficiently small (so that the leading term in the last line of the above display is positive) and hence conclude that $(\eta_n - \tau_n) = \mathrm{O}_p(n^{-1/5})$. A similar approach proves the non-decreasing case. □

**Lemma 5.4.** *Let $\xi_n$ be a sequence converging to $x_0$. Then, for any $\varepsilon > 0$, there exists an $M > 1$ and a $c > 0$ such that, with probability greater than $1 - \varepsilon$, we have that there exist change points $\tau_n < \xi_n < \eta_n$ of $\widehat{h}_n$ such that*

$$\inf_{t \in [\tau_n, \eta_n]} |\widehat{h}_n(t) - h_0(t)| < cn^{-2/5}$$

*for all $n$ sufficiently large.*

**Proof.** Fix $\varepsilon > 0$. From Lemma 5.3, it follows that there exist touchpoints $\eta_n$ and $\tau_n$, and an $M > 1$ such that $\xi_n - Mn^{-1/5} \leq \tau_n \leq \xi_n - n^{-1/5} \leq \xi_n + n^{-1/5} \leq \eta_n \leq \xi_n + Mn^{-1/5}$.

Fix $c > 0$ and consider the event

$$\inf_{t \in [\tau_n, \eta_n]} |\widehat{h}_n(t) - h_0(t)| \geq cn^{-2/5}. \tag{5.6}$$

First, assume that $\widehat{h}_n$ is non-increasing on $[\tau_n, \eta_n]$. On this set, we have that

$$\left| \int_{\tau_n}^{\eta_n} (\eta_n - t) \frac{\widehat{h}_n(t) - h_0(t)}{\widehat{h}_n(t)} S_0(t)\, dt \right|$$
$$\geq Bcn^{-2/5}(\eta_n - \tau_n)^2 \geq Bcn^{-4/5},$$

where $B$ is some constant depending on $x_0$. Using the definitions in (5.4), as well as the equality in condition (2.2) with (2.7), it follows that

$$0 = \widehat{\mathcal{H}}_{n,\downarrow}(\eta_n) - \int_0^{\eta_n} \mathbb{A}_{n,\downarrow}(t)\, dt - \widehat{\mathcal{H}}_{n,\downarrow}(\tau_n) + \int_0^{\tau_n} \mathbb{A}_{n,\downarrow}(t)\, dt$$
$$- (\widehat{\mathcal{H}}'_{n,\downarrow}(\tau_n) - \mathbb{A}_{n,\downarrow}(\tau_n))(\eta_n - \tau_n)$$
$$= \int_{\tau_n}^{\eta_n} \frac{\eta_n - t}{\widehat{h}_n(t)}\, d\mathbb{F}_n(t) - \int_{\tau_n}^{\eta_n} (\eta_n - t)\mathbb{S}_n(t)\, dt$$
$$= \int_{\tau_n}^{\eta_n} (\eta_n - t) \frac{\widehat{h}_n(t) - h_0(t)}{\widehat{h}_n(t)} S_0(t)\, dt + \int_{\tau_n}^{\eta_n} \frac{\eta_n - t}{\widehat{h}_n(t)}\, d\breve{\mathbb{F}}_n(t) - \int_{\tau_n}^{\eta_n} (\eta_n - t)\breve{\mathbb{S}}_n(t)\, dt,$$

where $\breve{\mathbb{F}}_n(t) = \mathbb{F}_n(t) - F_0(t)$ and $\breve{\mathbb{S}}_n(t) = \mathbb{S}_n(t) - S_0(t)$. By the assumption on $h_0$ and $x_0$, and arguments similar to those used for Lemmas 5.1 and 5.2, we can show that

$$\int_{\tau_n}^{\eta_n} (\eta_n - t) \frac{\widehat{h}_n(t) - h_0(t)}{\widehat{h}_n(t)} S_0(t)\, dt = O_p(n^{-4/5}),$$

which is a contradiction to (5.6) if $c$ is chosen large enough. A similar argument completes the proof for the non-decreasing case. □

The next proposition is the key to proving tightness in the next section. The results follow from the previous lemmas and make extensive use of the underlying convexity.

**Proposition 5.5.** *Under the assumptions of this section, we have that, for each $M > 0$,*

$$\sup_{|t| \leq M} |\widehat{h}_n(x_0 + n^{-1/5}t) - h_0(x_0) - n^{-1/5}t h_0'(x_0)| = O_p(n^{-2/5}), \tag{5.7}$$

$$\sup_{|t| \leq M} |\widehat{h}_n'(x_0 + n^{-1/5}t) - h_0'(x_0)| = O_p(n^{-1/5}). \tag{5.8}$$

**Proof.** For $M, \varepsilon > 0$ fixed, define $\eta_{n,1}$ to be the first point of touch after $x_0 + Mn^{-1/5}$ and $\eta_{n,i}$ to be the first point of touch after $\eta_{n,i-1} + n^{-1/5}$, $i = 2, 3$. Define the points $\eta_{n,-i}$ for $i = 1, 2, 3$ similarly, but working to the left of $x_0$. By Lemma 5.4, there exist points $\xi_{n,i} \in (\eta_{n,i}, \eta_{n,i+1})$, $i = 1, 2, -2, -3$, and a constant $c > 0$ such that with probability at least $1 - \varepsilon$, we have that $|\widehat{h}_n(\xi_{n,i}) - h_0(\xi_{n,i})| \leq cn^{-2/5}$.

As $\widehat{h}_n$ is convex, it follows that for any $t \in [x_0 - Mn^{-1/5}, x_0 + Mn^{-1/5}]$,

$$\widehat{h}_n'(t) \leq \widehat{h}_n'(\xi_{n,1}) \leq \frac{\widehat{h}_n(\xi_{n,2}) - \widehat{h}_n(\xi_{n,1})}{\xi_{n,2} - \xi_{n,1}}$$

$$\leq \frac{h_0(\xi_{n,2}) - h_0(\xi_{n,1}) + 2cn^{-2/5}}{\xi_{n,2} - \xi_{n,1}}$$

$$\leq h_0'(\xi_{n,2}) + 2cn^{-1/5}$$

since $\xi_{n,2} - \xi_{n,1} \geq n^{-1/5}$, where $\widehat{h}_n'(t)$ denotes the right derivative at $t$. Because of the continuity of $h_0''(\cdot)$ near $x_0$, we may replace $h_0'(\xi_{n,2})$ with $h_0'(x_0) + \tilde{c}n^{-1/5}$ for some new constant $\tilde{c}$. The result follows. A similar argument shows the lower bound.

We now consider (5.7). By Lemma 5.3, there exists a constant $K > M$ such that there exist two touchpoints in $[x_0 + Mn^{-1/5}, x_0 + Kn^{-1/5}]$, $n^{-1/5}$ apart with probability $1 - \varepsilon$. The same is the case in the interval $[x_0 - Mn^{-1/5}, x_0 - Kn^{-1/5}]$. From Lemma 5.4, it follows that there exist points $\xi_{n,1} \in [x_0 + Mn^{-1/5}, x_0 + Kn^{-1/5}]$ and $\xi_{n,2} \in [x_0 - Mn^{-1/5}, x_0 - Kn^{-1/5}]$ such that $|\widehat{h}_n(\xi_{n,i}) - h_0(\xi_{n,i})| \leq cn^{-2/5}$, for $i = 1, 2$, with probability at least $1 - \varepsilon$ and sufficiently large $n$. Lastly, we have already shown that there exists a $c'$ such that with probability at least $1 - \varepsilon$,

$$\sup_{t \in [x_0 - Kn^{-1/5}, x_0 + Kn^{-1/5}]} |\widehat{h}_n'(t) - h_0'(x_0)| \leq c'n^{-1/5}.$$

Therefore, with probability at least $1 - 3\varepsilon$, we have that for any $t \in [x_0 - Mn^{-1/5}, x_0 + Mn^{-1/5}]$ and sufficiently large $n$,

$$\widehat{h}_n(t) \geq \widehat{h}_n(\xi_{n,1}) + \widehat{h}_n'(\xi_{n,1})(t - \xi_{n,1})$$

$$\geq h_0(\xi_{n,1}) - cn^{-2/5} + (h_0'(x_0) - c'n^{-1/5})(t - \xi_{n,1})$$

$$= h_0(x_0) + h'(x_0)(t - x_0) + \tfrac{1}{2}h''(x_0^*)(\xi_{n,1} - x_0)^2 - cn^{-2/5} - c'n^{-1/5}(t - \xi_{n,1})$$

$$\geq h_0(x_0) + h'(x_0)(t - x_0) - Bn^{-2/5}$$

for some constant $B > 0$. A similar argument proves the other direction. □

### 5.3. Limit distribution theory at a fixed point

From Lemma 5.3, we know what rescaling is necessary to pick up a meaningful limit. The idea of the proof is now to write carefully a local version of the characterization of the MLE, Lemma 2.2, and to show that in the limit, these become the characterization of the invelope (Definition 1.1). The invelope $\mathcal{I}(\cdot)$ is described in terms of the "driving" process $Y(\cdot)$. Our goal will then be to identify the two processes, one which converges to the invelope and another which converges to the driving process $Y$.

Note that at $x_0$ (where $h''(x_0) > 0$), we have three possibilities:

1. $h_0'(x_0) > 0$: By continuity, $h_0'(x) > 0$ in a neighborhood of $x_0$. It follows from the consistency of the MLE derivatives that $\widehat{h}_n' > 0$ for sufficiently large $n$ and hence all touchpoints to be considered are of the "increasing" kind.
2. $h_0'(x_0) < 0$: By the same argument, all touchpoints are decreasing.
3. $h_0'(x_0) = 0$: Since $h(x_0) > 0$, by Corollary 2.5, there is always at least one touchpoint which satisfies both the non-increasing and non-decreasing properties. The limiting process may then be "stitched" together in an appropriate manner.

Therefore, it will be sufficient to prove the asymptotic results for both types of touchpoints: non-increasing and non-decreasing. For the sake of brevity, we outline the argument only for the non-increasing setting.

For any interval $[a,b] \subset \mathbb{R}$, let $D[a,b]$ denote the space of cadlag functions from $[a,b]$ into $\mathbb{R}$ endowed with the Skorohod topology and $C[a,b]$ the space of continuous functions endowed with the uniform topology.

*Driving process.*

Define

$$\mathbb{B}_n(t) \equiv \sqrt{n}(\mathbb{H}_n(t) - H_0(t)), \tag{5.9}$$

where $\mathbb{H}_n$ is the empirical cumulative hazard function, defined by $d\mathbb{H}_n(u) = (1 - \mathbb{F}_n(t-))^{-1} d\mathbb{F}_n(t)$. From [28], Chapter 7, Theorem 7.4.1, page 307, we know that for $t \in (0, T_0)$ with $T_0 \equiv T_0(F_0) \equiv \inf\{x : F(x) = 1\}$, $\mathbb{B}_n(t) \Rightarrow B(C(t))$ in $D[0, M]$ for $M < T_0$, where $B$ denotes a standard Brownian motion on $[0, \infty)$ and $C(t) = F_0(t)/S_0(t)$. Let $x_n(t) = x_0 + n^{-1/5}t$ and define

$$\widetilde{\mathbb{Y}}_n^{\mathrm{loc}}(t) \equiv n^{4/5} \int_{x_0}^{x_n(t)} \left\{ \mathbb{H}_n(v) - \mathbb{H}_n(x_0) - \int_{x_0}^{v} (h_0(x_0) + (u - x_0)h_0'(x_0)) \, du \right\} dv. \tag{5.10}$$

It is not difficult to show that

$$\widetilde{\mathbb{Y}}_n^{\mathrm{loc}}(t) \Rightarrow \sqrt{C'(x_0)} \int_0^t W(s) \, ds + \frac{1}{24} h_0''(x_0) t^4,$$

$$(\widetilde{\mathbb{Y}}_n^{\mathrm{loc}})'(t) \Rightarrow \sqrt{C'(x_0)}W(t) + \frac{1}{3!}h_0''(x_0)t^3$$

in $D[-M,M]$ for each fixed $0 < M < \infty$, where $W$ is a two-sided Brownian motion process starting at 0 and $C'(t) = h_0(t)/S_0(t)$. Next, define

$$\widehat{\mathbb{Y}}_{n,\downarrow}^{\mathrm{loc}}(t) = n^{4/5}\frac{h_0(x_0)}{S_0(x_0)}\int_{x_0}^{x_n(t)}\int_{x_0}^{v}\left\{\frac{h_0(u) - h_0(x_0) - (u - x_0)h_0'(x_0)}{\widehat{h}_n(u)}\right\}\mathbb{S}_n(u)\,\mathrm{d}u\,\mathrm{d}v$$

$$+ n^{4/5}\frac{h_0(x_0)}{S_0(x_0)}\int_{x_0}^{x_0+n^{-1/5}t}\int_{x_0}^{v}\frac{\mathbb{S}_n(u)}{\widehat{h}_n(u)}\,\mathrm{d}\{\mathbb{H}_n^*(u) - H_0(u)\}\,\mathrm{d}v,$$

where $\mathrm{d}\mathbb{H}_n^*(u) = \frac{\mathbb{S}_n(u^-)}{\mathbb{S}_n(u)}\,\mathrm{d}\mathbb{H}_n(u)$. The derivative $(\widehat{\mathbb{Y}}_{n,\downarrow}^{\mathrm{loc}})'(t)$ is not difficult to calculate. By consistency of $\widehat{h}_n$ and since $\sup_t|\mathbb{S}_n(t) - S_0(t)| \to 0$ a.s., for any $M > 0$,

$$\lim_n \sup_{|t| \leq M}|\widehat{\mathbb{Y}}_{n,\downarrow}^{\mathrm{loc}}(t) - \widetilde{\mathbb{Y}}_n^{\mathrm{loc}}(t)| = \lim_n \sup_{|t| \leq M}|(\widehat{\mathbb{Y}}_{n,\downarrow}^{\mathrm{loc}})'(t) - (\widetilde{\mathbb{Y}}_n^{\mathrm{loc}})'(t)| = 0 \qquad \text{a.s.} \quad (5.11)$$

$\widehat{\mathbb{Y}}_{n,\downarrow}^{\mathrm{loc}}$ is our driving process.

*Invelope process.*

Recall definitions (5.4). Our initial candidate for the invelope is defined as

$$\widehat{\mathbb{I}}_{n,\downarrow}^{\,\mathrm{loc}}(t) = n^{4/5}\frac{h_0(x_0)}{S_0(x_0)}\int_{x_0}^{x_n(t)}\int_{x_0}^{v}\left\{\frac{\widehat{h}_n(u) - h_0(x_0) - (u - x_0)h_0'(x_0)}{\widehat{h}_n(u)}\right\}\mathbb{S}_n(u)\,\mathrm{d}u\,\mathrm{d}v$$

$$+ \widehat{A}_{n,\downarrow}t + \widehat{B}_{n,\downarrow},$$

where

$$\widehat{A}_{n,\downarrow} = -n^{3/5}\frac{h_0(x_0)}{S_0(x_0)}\{\widehat{\mathcal{H}}_{n,\downarrow}'(x_0) - \mathbb{A}_{n,\downarrow}(x_0)\} \quad \text{and}$$

$$\widehat{B}_{n,\downarrow} = -n^{4/5}\frac{h_0(x_0)}{S_0(x_0)}\left\{\widehat{\mathcal{H}}_{n,\downarrow}(x_0) - \int_0^{x_0}\mathbb{A}_{n,\downarrow}(v)\,\mathrm{d}v\right\}.$$

Notice that because of the presence of $\mathbb{S}_n(v)$ in its definition, $\widehat{\mathbb{I}}_{n,\downarrow}^{\,\mathrm{loc}}(t)$ is not three times differentiable. We therefore define

$$\widehat{\mathbb{I}}_{n,\downarrow}^{\,*,\mathrm{loc}}(t) = n^{4/5}\frac{h_0(x_0)}{S_0(x_0)}\int_{x_0}^{x_n(t)}\int_{x_0}^{v}\left\{\frac{\widehat{h}_n(u) - h_0(x_0) - (u - x_0)h_0'(x_0)}{\widehat{h}_n(u)}\right\}S_0(u)\,\mathrm{d}u\,\mathrm{d}v$$

$$+ \widehat{A}_{n,\downarrow}t + \widehat{B}_{n,\downarrow}.$$

From Proposition 5.5, we have that for any $M > 0$,

$$\lim_n \sup_{|t| \leq M}|\widehat{\mathbb{I}}_{n,\downarrow}^{\,\mathrm{loc}}(t) - \widehat{\mathbb{I}}_{n,\downarrow}^{\,*,\mathrm{loc}}(t)| = 0. \tag{5.12}$$

The derivatives of $\widehat{\mathbb{I}}_{n,\downarrow}^{*,\mathrm{loc}}$ will describe the limiting behavior of our estimators. First, though, we must show that this process converges to the invelope. To do this, define the vector

$$\widehat{\mathbb{Z}}_n(t) = (\widehat{\mathbb{Y}}_{n,\downarrow}^{\mathrm{loc}}(t), (\widehat{\mathbb{Y}}_{n,\downarrow}^{\mathrm{loc}})'(t), \widehat{\mathbb{I}}_{n,\downarrow}^{*,\mathrm{loc}}(t), (\widehat{\mathbb{I}}_{n,\downarrow}^{*,\mathrm{loc}})'(t), (\widehat{\mathbb{I}}_{n,\downarrow}^{*,\mathrm{loc}})''(t), (\widehat{\mathbb{I}}_{n,\downarrow}^{*,\mathrm{loc}})'''(t)) \quad (5.13)$$

and fix $M > 0$. We will show that $\widehat{\mathbb{Z}}_n$ is tight in the product space

$$E[-M,M] \equiv C[-M,M] \times D[-M,M] \times C[-M,M]^3 \times D[-M,M].$$

This will be done last. We first assume that $\widehat{\mathbb{Z}}_n$ has a weak limit and identify its unique limit. The two arguments together prove that $\widehat{\mathbb{Z}}_n$, and hence $\widehat{\mathbb{I}}_{n,\downarrow}^{*,\mathrm{loc}}$, have the appropriate limiting distribution.

*Identifying the limit.*

It is sufficient to show that $\widehat{\mathbb{I}}_{n,\downarrow}^{*,\mathrm{loc}}$ satisfies (1.2)–(1.4) in the limit.

For condition (1.2), calculate

$$\widehat{\mathbb{I}}_{n,\downarrow}^{\mathrm{loc}}(t) - \widehat{\mathbb{Y}}_{n,\downarrow}^{\mathrm{loc}}(t) = n^{4/5} \frac{h_0(x_0)}{S_0(x_0)} \left\{ \int_0^{x_0 + n^{-1/5}t} \mathbb{A}_{n,\downarrow}(v)\,\mathrm{d}v - \widehat{\mathcal{H}}_{n,\downarrow}(x_0 + n^{-1/5}t) \right\} \geq 0, (5.14)$$

with equality at the (non-increasing) touchpoints of $\widehat{h}_n$, using (2.2). By (5.12), it follows that $\widehat{\mathbb{I}}_{n,\downarrow}^{*,\mathrm{loc}}$ satisfies (1.2) in the limit.

Next, the derivatives of $\widehat{\mathbb{I}}_{n,\downarrow}^{*,\mathrm{loc}}(t)$ are calculated as follows:

$$(\widehat{\mathbb{I}}_{n,\downarrow}^{*,\mathrm{loc}})'(t) = n^{3/5} \frac{h_0(x_0)}{S_0(x_0)} \int_{x_0}^{x_n(t)} \left\{ \frac{\widehat{h}_n(u) - h_0(x_0) - (u - x_0)h_0'(x_0)}{\widehat{h}_n(u)} \right\} S_0(u)\,\mathrm{d}u + \widehat{A}_{n,\downarrow},$$

$$(\widehat{\mathbb{I}}_{n,\downarrow}^{*,\mathrm{loc}})''(t) = n^{2/5} \frac{h_0(x_0)}{S_0(x_0)} \left\{ \frac{\widehat{h}_n(x_n(t)) - h_0(x_0) - n^{-1/5}th_0'(x_0)}{\widehat{h}_n(x_n(t))} \right\} S_0(x_0 + n^{-1/5}t).$$

Due to Theorem 3.1 and Proposition 5.5, we have that

$$\lim_n \sup_{|t| \leq M} |(\widehat{\mathbb{I}}_{n,\downarrow}^{*,\mathrm{loc}})''(t) - n^{2/5}[\widehat{h}_n(x_0 + n^{-1/5}t) - h_0(x_0) - n^{-1/5}th_0'(x_0)]| = 0, (5.15)$$

where $n^{2/5}[\widehat{h}_n(x_0 + n^{-1/5}t) - h_0(x_0) - n^{-1/5}th_0'(x_0)]$ is convex, and hence the limit of $(\widehat{\mathbb{I}}_{n,\downarrow}^{*,\mathrm{loc}})''(t)$ will be convex. Thus, (1.3) is satisfied in the limit.

Let $B_n(t) = (h_0(x_0)/S_0(x_0)) \times (S_0(t)/\widehat{h}_n(t))$. We may then write

$$(\widehat{\mathbb{I}}_{n,\downarrow}^{*,\mathrm{loc}})'''(t) = n^{1/5}[\widehat{h}_n'(x_n(t)) - h_0'(x_0)]B_n(x_n(t))$$
$$+ n^{1/5}[\widehat{h}_n(x_n(t)) - h_0(x_0) - n^{-1/5}th_0'(x_0)] \times B_n'(x_n(t)).$$

Notice that $\sup_{|t|\leq M}|1-B_n(x_0+n^{-1/5}t)| \to_{\text{a.s.}} 0$, with $\lim_n B'_n(x_0+n^{-1/5}t)$ bounded. Therefore, from Proposition 5.5, it follows that

$$\lim_n \sup_{|t|\leq M} |(\widehat{\mathbb{I}}_{n,\downarrow}^{*,\text{loc}})'''(t) - n^{1/5}[\widehat{h}'_n(x_0+n^{-1/5}t) - h'_0(x_0)]| = 0, \qquad (5.16)$$

where $n^{1/5}[\widehat{h}'_n(x_0+n^{-1/5}t) - h'_0(x_0)]$ is piecewise constant, with jumps at the touchpoints of $\widehat{h}_n$. By consistency of $\widehat{h}_n$, we have

$$\begin{aligned}
d(\widehat{\mathbb{I}}_{n,\downarrow}^{*,\text{loc}})'''(t) &= B_n(x_0+n^{-1/5}t)\,d\hat{g}_n(t) \\
&\quad + 2[\widehat{h}'_n(x_0+n^{-1/5}t) - h'_0(x_0)]B'_n(x_0+n^{-1/5}t)\,dt \\
&\quad + n^{1/5}[\widehat{h}_n(x_0+n^{-1/5}t) - h_0(x_0) - n^{-1/5}t h'_0(x_0)]\,dB'_n(x_0+n^{-1/5}t) \\
&= \{B_n(x_0+n^{-1/5}t) + O_p^*(n^{-2/5})\}\,d\hat{g}_n(t) + O_p^*(n^{-1/5})\,dt,
\end{aligned}$$

where $\hat{g}_n(t) = n^{1/5}[\widehat{h}'_n(x_0+n^{-1/5}t) - h'_0(x_0)]$. We say that a process $X_n(t)$ is $O_p^*(1)$ if $\sup_{|t|\leq M}|X_n(t)|$ is $O_p(1)$.

Next, fix a $c > 0$. Since $\widehat{h}_n$ is piecewise linear, it follows that $d\hat{g}_n$ puts mass only at the locations of touchpoints of $\widehat{h}_n$. However, at these locations, by (5.14), the process $\widehat{\mathbb{I}}_{n,\downarrow}^{\text{loc}}(t) - \widehat{\mathbb{Y}}_{n,\downarrow}^{\text{loc}}(t)$ is equal to zero. It follows that

$$\int_{-c}^c (\widehat{\mathbb{I}}_{n,\downarrow}^{\text{loc}}(t) - \widehat{\mathbb{Y}}_{n,\downarrow}^{\text{loc}}(t))\,d\hat{g}_n(t) = 0.$$

Hence,

$$\begin{aligned}
\int_{-c}^c (\widehat{\mathbb{I}}_{n,\downarrow}^{*,\text{loc}}(t) - \widehat{\mathbb{Y}}_{n,\downarrow}^{\text{loc}}(t))\,d(\widehat{\mathbb{I}}_{n,\downarrow}^{*,\text{loc}})'''(t) &= \int_{-c}^c (\widehat{\mathbb{I}}_{n,\downarrow}^{*,\text{loc}}(t) - \widehat{\mathbb{Y}}_{n,\downarrow}^{\text{loc}}(t))\,d[(\widehat{\mathbb{I}}_{n,\downarrow}^{*,\text{loc}})''' - \hat{g}_n](t) \\
&\quad + \int_{-c}^c (\widehat{\mathbb{I}}_{n,\downarrow}^{*,\text{loc}}(t) - \widehat{\mathbb{I}}_{n,\downarrow}^{\text{loc}}(t))\,d\hat{g}_n(t) = o_p(1),
\end{aligned}$$

using Proposition 5.5, (5.12) and the fact that $\hat{g}_n$ is increasing.

It remains to show that (1.4) is maintained under limits. This follows from the continuous mapping theorem since for any element $z = \{z_1, z_2, z_3, z_4, z_5, z_6\} \in E[-M, M]$,

$$\psi(z) = \int_{-M}^M (z_3 - z_1)\,dz_6$$

is continuous in $z$ for $z_6$ increasing. We have thus shown that $\widehat{\mathbb{I}}_{n,\downarrow}^{*,\text{loc}}(t)$ satisfies the invelope conditions (1.2)–(1.4) asymptotically. This shows that the only possible limit of $\widehat{\mathbb{I}}_{n,\downarrow}^{*,\text{loc}}(t)$ is the process $\mathcal{I}$.

*Tightness.*

We already know that $\widehat{\mathbb{Y}}_{n,\downarrow}^{\mathrm{loc}}(t)$ and $(\widehat{\mathbb{Y}}_{n,\downarrow}^{\mathrm{loc}})'(t)$ are tight in $C[-M,M]$ and $D[-M,M]$, respectively. To address tightness of the invelope processes, note that bounded and increasing functions are compact in $D[-M,M]$ and that bounded continuous functions with uniformly bounded derivatives are compact in $C[-M,M]$. These two facts allow us to address only stochastic boundedness of $(\widehat{\mathbb{I}}_{n,\downarrow}^{*,\mathrm{loc}})^{(i)}(t)$, $i=0,\ldots,3$, to obtain tightness. Thus, Proposition 5.5, along with (5.15) and (5.16), says that $(\widehat{\mathbb{I}}_{n,\downarrow}^{*,\mathrm{loc}})''(t)$ and $(\widehat{\mathbb{I}}_{n,\downarrow}^{*,\mathrm{loc}})'''(t)$ are tight. It remains to argue the same for $(\widehat{\mathbb{I}}_{n,\downarrow}^{*,\mathrm{loc}})'(t)$ and $\widehat{\mathbb{I}}_{n,\downarrow}^{*,\mathrm{loc}}(t)$. However, this will follow by Proposition 5.5, and (5.12), if we can show that both $\widehat{A}_{n,\downarrow}$ and $\widehat{A}_{n,\downarrow}t+\widehat{B}_{n,\downarrow}$ are tight.

Let $\tau_n$ be the largest touchpoint smaller than $x_0$. By (2.7), and after careful calculations, we have

$$-\frac{S_0(x_0)}{h_0(x_0)}\widehat{A}_{n,\downarrow} = -n^{3/5}\bigg\{\int_{\tau_n}^{x_0}\frac{\widehat{h}_n(u)-h_0(x_0)-h_0'(x_0)(u-x_0)}{\widehat{h}_n(u)}S_0(u)\,\mathrm{d}u\bigg\}$$
$$+n^{3/5}\bigg\{\int_{\tau_n}^{x_0}\frac{h_0(u)-h_0(x_0)-h_0'(x_0)(u-x_0)}{\widehat{h}_n(u)}S_0(u)\,\mathrm{d}u\bigg\}$$
$$+n^{3/5}\int_{\tau_n}^{x_0}\frac{1}{\widehat{h}_n(u)}\,\mathrm{d}\{\mathbb{F}_n-F_0\}(u)+n^{3/5}\int_{\tau_n}^{x_0}\mathbb{S}_n(u)-S_0(u)\,\mathrm{d}u.$$

By Proposition 5.5, Lemma 5.3 and Theorem 3.1, the first two terms are tight in $C[-M,M]$. Arguments similar to those used in the proof of Lemma 5.1, along with Lemma 5.3, may be used to handle the remaining terms. Since $\widehat{\mathcal{H}}_{n,\downarrow}(\tau_n)=\int_0^{\tau_n}\mathbb{A}_{n,\downarrow}(v)\,\mathrm{d}v$ by Lemma 2.2, it follows that $\widehat{B}_{n,\downarrow}$ is tight in $C[-M,M]$, which, in turn, implies that $\widehat{\mathbb{Z}}_n$ is tight in the space $E[-M,M]$.

*From the invelope to Theorem 1.2.*

By (5.15) and (5.16), the limiting behavior of $n^{2/5}(\widehat{h}_n(x_0)-h_0(x_0))$ and $n^{1/5}(\widehat{h}_n'(x_0)-h_0'(x_0))$ is the same as that of the second and third derivatives of $\widehat{\mathbb{I}}_{n,\downarrow}^{*,\mathrm{loc}}$, which converge to the invelope of $\lim_n \widehat{\mathbb{Y}}_{n,\downarrow}^{\mathrm{loc}}$. Define $k_1,k_2$ by

$$\lim_n \widehat{\mathbb{Y}}_{n,\downarrow}^{\mathrm{loc}}(t) = \sqrt{C'(x_0)}\int_0^t W(s)\,\mathrm{d}s + \frac{1}{24}h_0''(x_0)t^4 \equiv k_1\int_0^t W(s)\,\mathrm{d}s + k_2 t^4.$$

For any $a,b>0$, $bY(at)\stackrel{d}{=}a^{3/2}b\int_0^t W(s)\,\mathrm{d}s+a^4bt^4$. Therefore, choose $a,b$ so that $a^4b=k_2$ and $a^{3/2}b=k_1$. It follows that

$$\widehat{\mathbb{Y}}_{n,\downarrow}^{\mathrm{loc}}(t)\Rightarrow bY(at).$$

Applying this rescaling to all processes shows that

$$(\widehat{\mathbb{I}}_{n,\downarrow}^{*,\mathrm{loc}})''(0)\Rightarrow ba^2\mathcal{I}''(0)\quad\text{and}\quad(\widehat{\mathbb{I}}_{n,\downarrow}^{*,\mathrm{loc}})'''(0)\Rightarrow ba^3\mathcal{I}'''(0).$$

It is now straightforward to calculate the correct constants, $c_1$ and $c_2$, of Theorem 1.2.

## Acknowledgements

The research of Hanna Jankowski was supported by the NSERC; the research of Jon A. Wellner was supported in part by NSF Grant DMS-0503822 and NI-AID Grant 2R01 AI291968-04.